\documentclass{amsart}
\usepackage{amssymb,latexsym}
\numberwithin{equation}{section}
\newtheorem{theorem}{Theorem}[section]
\newtheorem{proposition}[theorem]{Proposition}

\newtheorem{lemma}[theorem]{Lemma}

\newtheorem{example}[theorem]{Example}

\begin{document}

\pagenumbering{arabic}
\pagestyle{headings}
\def\sof{\hfill\rule{2mm}{2mm}}
\def\ls{\leq}
\def\gs{\geq}
\def\SS{\mathcal S}
\def\qq{{\bold q}}
\def\txx{{\frac1{2\sqrt{x}}}}
\def\mn{\mbox{-}}
\def\dd{\mbox{-}}

\def\axbc{\ensuremath{1{\dd}23}}
\def\axcb{\ensuremath{1{\dd}32}}
\def\bxac{\ensuremath{2{\dd}13}}
\def\bxca{\ensuremath{2{\dd}31}}
\def\cxab{\ensuremath{3{\dd}12}}
\def\cxba{\ensuremath{3{\dd}21}}

\def\abxc{\ensuremath{12{\dd}3}}
\def\acxb{\ensuremath{13{\dd}2}}
\def\baxc{\ensuremath{21{\dd}3}}
\def\bcxa{\ensuremath{23{\dd}1}}
\def\caxb{\ensuremath{31{\dd}2}}
\def\cbxa{\ensuremath{32{\dd}1}}

\def\ls{\leqslant}
\def\gs{\geqslant}
\def\flm{f_{l,m}}
\def\plm{P_{l,m}}
\def\flma{f_{l,m}^a}
\def\Flma{F_{l,m}^a}

\title{ The centralizer of two numbers under the natural action of $S_k$ on
$[k]$, the maximal parabolic subgroup of $S_k$, and generalized
patterns }

\author[Toufik Mansour]{T. Mansour } \maketitle
\begin{center}
{\small LaBRI (UMR 5800), Universit\'e Bordeaux 1,
       351 cours de la Lib\'eration\\ 33405 Talence Cedex, France}\\[4pt]
{\tt toufik@labri.fr}
\end{center}
\markboth{\large Toufik Mansour} {\large Restricted certain
generalized patterns}
%
%===========================================================================
\section*{Abstract}
A natural generalization of single pattern avoidance is {\it
subset avoidance\/}. A complete study of subset avoidance for the
case $k=3$ is carried out in \cite{SS}. For $k>3$ situation
becomes more complicated, as the number of possible cases grows
rapidly. Recently, several authors have considered the case of
general $k$ when $T$ has some nice algebraic properties. Barcucci,
Del Lungo, Pergola, and Pinzani in \cite{BDPP} treated the case
when $T=T_1$ is the centralizer of $k-1$ and $k$ under the natural
action of $S_k$ on $[k]$. Mansour and Vainshtein in \cite{MVp}
treated the case when $T=T_2$ is maximal parabolic group of $S_k$.

Recently, Babson and Steingrimsson (see \cite{BS}) introduced
generalized permutations patterns that allow the requirement that
two adjacent letters in a pattern must be adjacent in the
permutation.

In this paper we present an analogue with generalization for the
case $T_1$ and for the case $T_2$ by using generalized patterns
instead of classical patterns.
%===========================================================================
\section{Introduction}

{\bf Classical patterns}. Let $[p]=\{1,\dots,p\}$ denote a totally
ordered alphabet on $p$ letters, and let
$\alpha=(\alpha_1,\dots,\alpha_m)\in [p_1]^m$,
$\beta=(\beta_1,\dots,\beta_m)\in [p_2]^m$. We say that $\alpha$
is {\it order-isomorphic\/} to $\beta$ if for all $1\ls i<j\ls m$
one has $\alpha_i<\alpha_j$ if and only if $\beta_i<\beta_j$. For
two permutations $\pi\in S_n$ and $\tau\in S_k$, an {\it
occurrence\/} of $\tau$ in $\pi$ is a subsequence $1\ls
i_1<i_2<\dots<i_k\ls n$ such that $(\pi_{i_1}, \dots,\pi_{i_k})$
is order-isomorphic to $\tau$; in such a context $\tau$ is usually
called the {\it pattern\/}. We say that $\pi$ {\it avoids\/}
$\tau$, or is $\tau$-{\it avoiding\/}, if there is no occurrence
of $\tau$ in $\pi$. Pattern avoidance proved to be a useful
language in a variety of seemingly unrelated problems, from stack
sorting \cite[Ch.~2.2.1]{Kn} to singularities of Schubert
varieties \cite{LS}. A natural generalization of single pattern
avoidance is {\it subset avoidance\/}; that is, we say that
$\pi\in S_n$ avoids a subset $T\subset S_k$ if $\pi$ avoids any
$\tau\in T$. A complete study of subset avoidance for the case
$k=3$ is carried out in \cite{SS}. For $k>3$ situation becomes
more complicated, as the number of possible cases grows rapidly.
Recently, several authors have considered the case of general $k$
when $T$ has some nice algebraic properties.

Adin and Roichman in \cite{AR} treated the case when $T$ is a
Kazhdan--Lusztig cell of $S_k$, or, equivalently, the Knuth
equivalence class (see \cite[vol.~2, Ch.~A1]{St}). Barcucci, Del
Lungo, Pergola, and Pinzani in \cite{BDPP} treated the case when
$T$ is the centralizer of $k-1$ and $k$ under the natural action
of $S_k$ on $[k]$.

\begin{theorem} {\rm (Barcucci, Del
Lungo, Pergola, and Pinzani ~\cite{BDPP})} \label{bd1} Let $T$ be
the centralizer of $k-1$ and $k$ under the natural action of $S_k$
on $[k]$. Then the ordinary generating function for the number of
$T$-avoiding permutations in $S_n$ is given by
$$(k-3)!x^{k-4}\frac{1-(k-1)x-\sqrt{1-2(k-1)x+(k-3)^2x^2}}{2}+\sum_{i=0}^{k-3}i!x^i.$$
\end{theorem}

In 2000, Kremer \cite{Kre} presented a generalization for the
above theorem, but later Mansour \cite[Section~2.3.3]{Mathesis}
found a counterexample to the main result in \cite{Kre}.\\

Recall that the {\it Laguerre polynomial\/} $L_n^\alpha(x)$ is
given by
$$
L_n^\alpha(x)=\frac 1{n!}e^xx^{-\alpha}\frac{d^n}{dx^n}
\left(e^{-x}x^{n+\alpha}\right),
$$
and the {\it rook polynomial\/} of the rectangular $s\times t$
board is given by
$$
R_{s,t}(x)=s!x^sL_s^{t-s}(-x^{-1})
$$
for $s\ls t$ and by $R_{s,t}(x)=R_{t,s}(x)$ otherwise (see
\cite[Ch.~7.4]{Ri}). In 2000, Mansour and Vainshtein \cite{MVp}
treated the case when we $T$ is a maximal parabolic subgroup of
$S_k$.

\begin{theorem} {\rm(Mansour and Vainshtein ~\cite{MVp})} \label{par}
Let $1\leq a\leq m+1$, $\lambda=\min\{l,m\}$, $\mu=\max\{l,m\}$,
and let ${P'}_{l,m}=\{\sigma\in S_{l+m}| a\leq \sigma_j\leq
a+l-1,\ j=1,2,\dots,l\}$, then
$$
\Flma(x)R_{l,m}(-x)=\sum_{r=0}^{\lambda-1}x^rr!
\sum_{j=0}^r(-1)^j\frac{\binom lj\binom mj}{\binom rj}+
(-1)^\lambda x^\lambda\lambda!\sum_{r=0}^{\mu-\lambda-1}x^rr!
\binom {\mu-r-1}\lambda,
$$
or, equivalently,
$$
\Flma(x)=\sum_{r=0}^{k-1}x^rr!-\frac{(-1)^\lambda x^\mu}
{\lambda!L_\lambda^{\mu-\lambda}(x^{-1})}\sum_{r=0}^{\lambda-1}
(k+r)!x^r\sum_{j=r+1}^\lambda(-1)^j\frac{\binom lj\binom
mj}{\binom {k+r}j},
$$
where $k=l+m=\lambda+\mu$, $\Flma(x)$ is the ordinary generating
function for the number of ${P'}_{l,m}$-avoiding permutations in
$S_n$.\\
\end{theorem}
%==============
{\bf Generalized patterns}. In \cite{BS} Babson and
Steingr\'{\i}msson introduced generalized permutation patterns
that allow the requirement that two adjacent letters in a pattern
must be adjacent in the permutation. For example, (generalized)
patterns are $123\mn4$ and $12\mn34$. An occurrence of $123\mn4$
in a permutation $\pi$ is a subword $\pi_i
\pi_{i+1}\pi_{i+2}\pi_j$ of $\pi$ such that
$\pi_i<\pi_{i+1}<\pi_{i+2}<\pi_j$, and an occurrence of $12\mn34$
is a subword $\pi_i\pi_{i+1}\pi_{j}\pi_{j+1}$ of $\pi$ such that
$\pi_i<\pi_{i+1}<\pi_{j}<\pi_{j+1}$.
%More generally, if $\sigma\in
%S_k$, then we define an occurrence of
%$\sigma_1\sigma_2\dots\sigma_d\mn\sigma_{d+1}\mn\dots\mn\sigma_k$
%in $\pi$ is a subword
%$\pi_i\pi_{i+1}\dots\pi_{i+d-1}\pi_{j_{d+1}}\pi_{j_{d+2}}\dots\pi_{j_k}$
%order-isomorphic to $\sigma$, where $i+d\leq
%j_{d+1}<j_{d+2}<\dots<j_k$.

Claesson in \cite{C} presented a complete solution for the number
of permutations avoiding any single (generalized) pattern of
length three with exactly one adjacent pair of letters (we do not
regard a dash as being a letter) as follows.

\begin{proposition}{\rm(Claesson ~\cite{C})}\label{c1}
  For all $n\geq 0$
  $$
  |S_n(\sigma)| =\left\{
  \begin{array}{ll}
    B_n\quad & \text{ if }\,
    \sigma\in\{
    \axbc, \cxba, \abxc, \cbxa,
    \axcb, \cxab, \baxc, \bcxa
    \},\\
    C_n &\text{ if }\,
    \sigma\in\{
    \bxac, \bxca, \acxb, \caxb
    \},
  \end{array}
  \right.
  $$
  where $B_n$ and $C_n$ are the $n$th Bell and Catalan numbers,
  respectively.
\end{proposition}
In addition, Claesson in \cite{C} gave certain results for the
number of permutations avoiding a pair generalized patterns of
three letters.
\begin{proposition}{\rm(Claesson ~\cite{C})}\label{claesson2}
  For all $n\geq 0$
  $$
  S_n(\axbc,\,\abxc) = B^*_n,\;\;\;
  S_n(\axbc,\,\axcb) = I_n,\;\;\text{ and }\;\,
  S_n(\axbc,\,\acxb) = M_n,
  $$
  where $B^*_n$ is the $n$th Bessel number \textup{(}\#
  non-overlapping partitions of $[n]$ \textup{(}see
  \cite{Fl80}\textup{)}\textup{)}, $I_n$ is the number of involutions
  in $S_n$, and $M_n$ is the $n$th Motzkin number.\\
\end{proposition}

Later, Claesson and Mansour ~\cite{CM1} gave the complete answer
for the number permutations avoiding a pair patterns of the form
$ab\mn c$ or $a\mn bc$ where $abc\in S_3$.\\

In this paper we present an analogue for Theorem~\ref{bd1} and for
Theorem~\ref{par} by using generalized patterns instead of
classical patterns as following.\\

Let $1\leq a<a+l\leq k$ and let us denote by $C_{a,l}^k$ the set
of all generalized patterns
$\sigma_1\sigma_2\mn\sigma3\mn\dots\mn\sigma_k$ such that
$(\sigma_1,\dots,\sigma_k)\in S_k$, $\sigma_1=a$, and
$\sigma_2=a+l$. Clearly, the set $C_{a,l}^k$ contains $(k-2)!$
generalized patterns ($C_{a,l}^k$ is an analogue for the case when
$T$ is the centralizer of $a$ and $a+l$ under the natural action
of $S_k$ on $[k]$; that is, $T=\{\pi\in S_k|
\pi_1=a,\,\pi_2=a+l\}$). For example, $C_{1,2}^4=\{13\mn2\mn4,
13\mn4\mn2\}$. In the present paper we find the number of
$C_{a,l}^k$-avoiding permutations in $S_n$ (see
Section~\ref{sec2}).

Let $P_{a,l}^k$ be all the generalized patterns of the form
$\sigma_1\dots\sigma_l\mn\sigma_{l+1}\mn\dots\mn\sigma_{k}$ such
that $(\sigma_1,\dots,\sigma_l)$ is a permutation of the numbers
$a,\dots,a+l-1$ and $(\sigma_{l+1},\dots,\sigma_{k})$ is a
permutation of the numbers $a+l,\dots,k,1,2,\dots,a-1$. Clearly,
$P_{a,l}^k$ contains $l!\cdot(k-l)!$ generalized patterns
($P_{1,l}^k$ is an analogue for the case when $T$ is a maximal
parabolic subgroup of $S_{k}$; that is, $T=\{\pi\in
S_k|(\pi_1,\dots,\pi_l)\in S_l\}$). For example, $P_{2,2}^4$
contains $4$ generalized patterns which are $23\mn1\mn4$,
$23\mn4\mn1$, $32\mn1\mn4$, and $32\mn4\mn1$. In the present paper
we find the number of $P_{a,l}^k$-avoiding permutations in $S_n$
(see Section~\ref{sec3}).
%===============================================================
\section{$C_{a,l}^k$-Avoiding}\label{sec2}

Let $c_{a,l}^k(n)$ be the number of $C_{a,l}^k$-avoiding
permutations in $S_n$. Our present aim is to count this number,
since that we introduce another notation. We denote by
$c_{a,l}^k(n; i_1,i_2,\ldots,i_j)$ the number of
$C_{a,l}^k$-avoiding permutations $\pi\in S_n$ such that
$\pi_1\pi_2\dots\pi_j=i_1i_2\dots i_j$. Now we introduce the two
quantities that play the crucial role in the proof of the main
theorem in this section.

\begin{lemma}\label{ls2}
For all $n\geq k$,
    $$c_{a,l}^k(n)=(k-1)c_{a,l}^k(n-1)+\sum_{j=a}^{n-k+a}c_{a,l}^k(n;j).$$
\end{lemma}
\begin{proof}
Let $\pi\in S_n$ where $\pi_1=i$, and let $n\geq j\geq n-k+a+1$ or
$1\leq j\leq a-1$. By definitions it is easy to see that every
occurrence of $\tau\in C_{a,l}^k$ in $\pi$ not contains $\pi_1$,
so $c_{a,l}^k(n;j)=c_{a,l}^k(n-1)$. On the other hand, by
definitions $c_{a,l}^k(n)=\sum_{j=1}^n c_{a,l}^k(n;j)$, hence the
lemma holds.
\end{proof}

\begin{lemma}\label{ls1}
Let $n\geq k$, and $a\leq j\leq n-k+a$. Then
$$c_{a,l}^k(n;j)=(k-l-1)c_{a,l}^k(n-2)+\sum_{i=a}^{j+l-2}c_{a,l}^k(n-1;i).$$
\end{lemma}
\begin{proof}
Let $a\leq j\leq n-k+a$; by definitions
$$c_{a,l}^k(n;j)=\sum_{i=1}^{j-1}c_{a,l}^k(n;j,i)+\sum_{i=j+1}^nc_{a,l}^k(n;j,i),\eqno(\ast)$$
and let us consider the possible values of $c_{a,l}^k(n;j,i)$:

\begin{itemize}
\item[$(i)$] Let $1\leq i\leq a-1$; since all the generalized patterns in
$C_{a,l}^k$ started by the segment $a(a+l)$ we get
$$c_{a,l}^k(n;j,i)=c_{a,l}^k(n-2);$$

\item[$(ii)$] Let $a\leq i\leq j-1$; similarly as $(i)$ with $a<a+l$ we
have
$$c_{a,l}^k(n;j,i)=c_{a,l}^k(n-1;i);$$

\item[$(iii)$] Let $j+1\leq i\leq j+l-1$; similarly as $(i)$ we
obtain
$$c_{a,l}^k(n;j,i)=c_{a,l}^k(n-1;i-1);$$

\item[$(iv)$] Let $j+l\leq i\leq n-k+a+l$; so there exist $a_1,\dots,a_{k-2}$
positions such that $(j,i,\pi_{a_1},\dots,\pi_{a_{k-2}})$ is
order-isomorphic to $\tau\in C_{a,l}^k$; hence
$$c_{a,l}^k(n;j,i)=0;$$

\item[$(v)$] Let $n-k+a+l+1\leq j\leq n$; similarly as $(i)$ with $a<a+l$ we get
$$c_{a,l}^k(n;j,i)=c_{a,l}^k(n-2).$$
\end{itemize}

Hence, if summing $c_{a,l}^k(n;j,i)$ over all $i$ with using
($\ast$) we get the desired result.
\end{proof}

Now we ready to obtain the main quantity that plays the crucial
role in the proof of the main theorem in this section.

\begin{proposition}\label{ml}
Let $n\geq k$, and $1\leq i\leq n-k+1$. Then
    $$c_{a,l}^k(n;n-k+a+1-i)=\sum_{j=0}^{\lfloor (i-1)/l\rfloor+1}(-1)^j\binom{i-(j-1)(l-1)}{j}c_{a,l}^k(n-1-j).$$
\end{proposition}
\begin{proof}
Using Lemma~\ref{ls2} for $i=1$ we get
$$c_{a,l}^k(n;n-k+a)=(k-l-1)c_{a,l}^k(n-2)+\sum_{j=a}^{n-k+b-2}c_{a,l}^k(n-1;j).$$
Using Lemma~\ref{ls1} with $c_{a,l}^k(n-1;j)=c_{a,l}^k(n-2)$ for
$1\leq j\leq a-1$ or $n-k+a+1\leq j\leq n$ (see Lemma~\ref{ls2}),
and with $a+l>a$ we have
$$c_{a,l}^k(n;n-k+a)=(k-l-1)c_{a,l}^k(n-2)+c_{a,l}^k(n-1)-(a-1+k-a-l+1)c_{a,l}^k(n-2),$$
equivalently
$$c_{a,l}^k(n;n-k+a)=c_{a,l}^k(n-1)-c_{a,l}^k(n-2).$$
Therefore the proposition holds for $i=1$.

Now, let $i\geq 2$; using Lemma~\ref{ls2} we obtain for all
$a+1\leq m\leq n-k+a$
$$c_{a,l}^k(n;m-1)=c_{a,l}^k(n;m)-c_{a,l}^k(n-1;m+l-2),$$
and by means of induction on $i$ with employing the familiar
identity
$$\binom{0}{s}+\binom{1}{s}+\binom{2}{s}+\cdots+\binom{t}{s}=\binom{t+1}{s+1},\eqno(\ast\ast)$$
we get the desired result.
\end{proof}

Now we are ready to prove the main result of this section.

\begin{theorem}\label{main11}
For all $n\geq k$
$$\begin{array}{l}
c_{a,l}^k(n)=(k-l-1)c_{a,l}^k(n-1)+\\
\qquad\qquad\qquad+\sum\limits_{j=0}^{\lfloor(n-k)/l\rfloor+1}(-1)^j\binom{n-k+2-(j-1)(l-1)}{j+1}c_{a,l}^k(n-1-j).
\end{array}$$
\end{theorem}
\begin{proof}
Lemma~\ref{ls1} and Proposition~\ref{ml} yield for all $n\geq k$
$$c_{a,l}^k(n)=(k-1)c_{a,l}^k(n-1)+\sum_{i=1}^{n-k+1}\sum_{j=0}^{\lfloor(i-1)/l\rfloor+1}(-1)^j\binom{i-(j-1)(l-1)}{j}c_{a,l}^k(n-1-j).$$
Again, employing the familiar identity ($\ast\ast$) we get the
desired result.
\end{proof}

\begin{example}{\rm(see Claesson and Mansour~\cite{CM1} for $k=3$)}
Theorem~\ref{main11} yields for $l=1$
   $$(k-2)c_{a,1}^k(n-1)=\sum_{j=0}^{n-k+2}(-1)^{n-k+2-j}\binom{n-k+2}{j}c_{a,1}^k(j+k-2).$$
Using \cite[Lem.~7]{CM2} we have
$$\begin{array}{l}
(k-2)uC_{a,1}^k(u)=(1+x)^{k-3}\left[C_{a,1}^k\left(\frac{u}{1+u}\right)-1\right]+\\
\qquad\qquad\qquad\qquad\qquad\sum\limits_{j=1}^{k-3}(j-1)!u^j(k-2-j(1+u)^{k-3-j}),
\end{array}$$ where $C_{a,1}^k(x)$ is the ordinary generating function for
the sequence $\{c_{a,1}^k(n)\}_{n\geq 0}$, hence by putting
$u=x/(1-x)$ we get
$$C_{a,1}^k(x)=1+(k-2)x(1-x)^{k-4}C_{a,1}^k\left(\frac{x}{1-x}\right)+\sum_{j=1}^{k-3}(j-1)!x^j(j-(k-2)(1-x)^{k-3-j}).$$
An infinite number of applications of this identity we have
$$\begin{array}{l}
C_{a,1}^k(x)=\sum\limits_{n\geq0}\frac{(k-2)^nx^n(1-nx)^{k-4}}{(1-x)(1-2x)\dots(1-(n-1)x)}+\\
\qquad\qquad\qquad+\sum\limits_{j=1}^{k-3}\sum\limits_{n\geq0}\frac{j!(k-2)^nx^{j+n}(1-nx)^{k-4-j}}{(1-x)(1-2x)\cdots(1-(n-1)x)}+\\
\qquad\qquad\qquad\qquad\qquad+\sum\limits_{j=0}^{k-4}\sum\limits_{n\geq0}\frac{j!(k-2)^{n+1}x^{j+n+1}(1-(n+1)x)^{k-4-j}}{(1-x)(1-2x)\cdots(1-nx)}.
\end{array}$$
An example, for $k=3$ we get
    $$C_{1,1}^3(x)=C_{2,1}^3(x)=\sum_{n\geq 0}\frac{x^n}{(1-x)(1-2x)\cdots(1-(n+1)x)}.$$
\end{example}

In view of Theorem~\ref{main11} we get the number of
$C_{1,l}^k$-avoiding permutations in $S_n$ is the same number of
$C_{a,l}^k$-avoiding permutations in $S_n$ where
$a=1,2,\dots,k-l$. Our proof in Theorem~\ref{main11} is an
analytical proof. It is a challenge to find a bijective proof for
this property.

To put Theorem~\ref{main11} in the form of ordinary generating
function let us define as follows. Let $l\geq 2$, we define for
all $n\geq0$,
$$F_n^l(x)=\sum_{m\geq0} (-1)^j\binom{n-(l-2)j}{j}x^j.$$
For example, $F_n^2(x)=(1-x)^n$.

From now we assume that $l\geq 2$. If multiplying the recurrence
in the statement of Theorem~\ref{main11} by $x^n$ and summing over
all $n\geq k$, then
$$\begin{array}{l}
(k-l-1)\sum\limits_{n\geq k}c_{a,l}^k(n-1)x^n=
%(k-2)!x^k-(k-1)!x^k(l+1-3x+\delta{l,2}x^2)+\\
\sum\limits_{n\geq k}\sum\limits_{j=0}^{n+1}
(-1)^j\binom{n-k+2-(j-1)(l-1)}{j+1}c_{a,l}^k(n-1-j)x^{n},
\end{array}$$ equivalently,

\begin{theorem}\label{main12}
Let $l\geq 2$, then
$$\begin{array}{l}
x(k-1-l)\sum\limits_{n\geq k-1}c_{a,l}^k(n)x^n=(k-1)!(F_{2l-1}^l(x)-1)x^{k-1}+\\
\qquad\qquad\qquad+(k-2)!(F_{2l-2}^l(x)-1+(2l-2)x)x^{k-2}+\sum\limits_{n\geq
k}c_{a,l}^k(n)F_{n-k+2l}^l(x)x^n.
\end{array}$$
\end{theorem}

For example (see Claesson ~\cite{C}), for $l=2$ and $k=3$
Theorem~\ref{main12} yields $c_{1,2}^3(n)$ is given by the $n$th
Catalan number; that is, the number of permutations in
$S_n(13\mn2)$ is $\frac{1}{n+1}\binom{2n}{n}$.
%========================================================================
\section{$P_{a,l}^k$-Avoiding}\label{sec3}
Let $p_{a,l}^k(n)$ be the number $P_{a,l}^k$-avoiding permutations
in $S_n$. Our present aim is to count this number, since that we
introduce another notation. We denote by
$p_{a,l}^k(n;i_1,\dots,i_m)$  the number $P_{a,l}^k$-avoiding
permutations $\pi\in S_n$ such that $\pi_1\dots\pi_m=i_1\dots
i_m$. Now we introduce the quantity that play the crucial role to
find $p_{a,l}^k(n)$.

\begin{lemma}\label{lp1}
Let $a\leq i_j\leq n-k+l+a-1$ for any $j=1,2,\dots,m$.\\

{\rm (i)} If $0\leq m\leq l-1$, then
$p_{a,l}^k(n;i_1,\dots,i_m,j)=p_{a,l}^k(n-m-1)$ where either
$1\leq j\leq a-1$, or $n-k+l+a\leq j\leq n$;\\

{\rm (ii)} If $m=l$, then $p_{a,l}^k(n;i_1,\dots,i_l)=0$;\\

{\rm (iii)} For all $i$, $p_{a,l}^k(n;\dots,i,\dots,i,\dots)=0$.
\end{lemma}
\begin{proof}
(iii) yields immediately by definitions. To verify (i) let $\pi\in
S_n$ a permutations such that $\pi_j=i_j$ for $j=1,2,\dots,m$, and
$\pi_{m+1}=j$. Let $\tau\in P_{a,l}^k$ and let
$(\pi_{a_1},\dots,\pi_{a_k})$ order-isomorphic to $\tau$. Since
$m\leq l-1$, and either $j\leq a-1$ or $j\geq n-k+l+a$ we get
$a_1>m+1$. Therefore, $p_{a,l}^k(n;i_1,\dots,i_m,j)$ is equal to
the number of $P_{a,l}^k$-avoiding permutations on the letters
$1,2,\dots,n$ without the letters $i_1,\dots,i_m,j$. Hence (i)
holds.

To verify (ii) let $\pi\in S_n$ a permutation such that
$\pi_j=i_j$ for all $j=1,2,\dots,l$, so there exist
$a_1,\dots,a_{k-l}$ positions such that
$(\pi_{i_1},\dots,\pi_{i_l},\pi_{a_1},\dots,\pi_{a_{k-l}})$ is
order-isomorphic to $\tau$ where $\tau\in P_{a,l}^k$. Hence
$p_{a,l}^k(n;i_1,\dots,i_l)=0$.
\end{proof}

An application of the above lemma we find the exponential
generating function for the sequence $\{p_{a,l}^k(n)\}_{n\geq 0}$
as follows.

\begin{theorem}\label{tp1}
Let $k,l,a\geq 1$ and $1\leq a+l-1\leq k$. Then
$$\sum_{n\geq 0}p_{a,l}^k(n)\frac{x^n}{n!}=\underbrace{\int\int\cdots\int}_{(k-l-1)-\mbox{times}}
e^{(k-l)(x/1+x^2/2+\dots+x^l/l)}dx\cdots dx.$$
\end{theorem}
\begin{proof}
By definitions we have
$$p_{a,l}^k(n)=\sum_{i_1=1}^n p_{a,l}^k(n;i_1).$$
Using Lemma~\ref{lp1}(i) for $m=0$ we get
$$p_{a,l}^k(n)=(k-l)p_{a,l}^k(n-1)+\sum_{i_1=a}^{n-k+l+a-1} p_{a,l}^k(n;i_1).$$
Now, let us assume for $m\geq 2$
    $$p_{a,l}^k(n)-(k-l)\sum_{j=0}^{m-1}j!\binom{n-k+l}{j}p_{a,l}^k(n-1-j)=\sum_{i_1,\dots,i_m=1}^{n-k+l+a-1} p_{a,l}^k(n;i_1,\dots,i_m),$$
so by definitions we have
$$\begin{array}{l}
p_{a,l}^k(n)-(k-l)\sum\limits_{j=0}^{m-1}j!\binom{n-k+l}{j}p_{a,l}^k(n-1-j)=\\
\qquad\qquad\qquad\qquad\qquad=\sum\limits_{i_1,\dots,i_m=1}^{n-k+l+a-1}\sum\limits_{i_{m+1}=1}^np_{a,l}^k(n;i_1,\dots,i_m,i_{m+1}),
\end{array}$$ equivalently
$$\begin{array}{l}
p_{a,l}^k(n)-(k-l)\sum\limits_{j=0}^{m-1}j!\binom{n-k+l}{j}p_{a,l}^k(n-1-j)=\\
\quad=\sum\limits_{i_1,\dots,i_m+1=1}^{n-k+l+a-1}p_{a,l}^k(n;i_1,\dots,i_m,i_{m+1})
+\sum\limits_{j=1}^{a-1}\sum\limits_{i_1,\dots,i_m=1}^{n-k+l+a-1}p_{a,l}^k(n;i_1,\dots,i_m,j)+\\
\qquad+\sum\limits_{j=n-k+l+a}^n\sum\limits_{i_1,\dots,i_m=1}^{n-k+l+a-1}p_{a,l}^k(n;i_1,\dots,i_m,j).
\end{array}$$
Using Lemma~\ref{lp1}(i)-(iii) we get
$$\begin{array}{l}
p_{a,l}^k(n)-(k-l)\sum\limits_{j=0}^{m-1}j!\binom{n-k+l}{j}p_{a,l}^k(n-1-j)=\\
\qquad=\sum\limits_{i_1,\dots,i_m+1=1}^{n-k+l+a-1}p_{a,l}^k(n;i_1,\dots,i_m,i_{m+1})+
(k-l)m!\binom{n-k+l}{m}p_{a,l}^k(n-1-m).
\end{array}$$ Hence, by the principle of induction with $m=l$ with Lemma ~\ref{lp1} we get
for all $n\geq 1$
    $$p_{a,l}^k(n)=(k-l)\sum_{j=0}^{l-1}j!\binom{n-k+l}{j}p_{a,l}^k(n-1-j).$$
Let $f(x)=\sum_{n\geq 0}p_{a,l}^k(n)x^n/n!$; if multiplying by
$x^n/(n-k+l)!$ and summing over all $n\geq 1$ we get
     $$\frac{{\rm d}^{k-l}}{{\rm d}x^{k-l}}f(x)=(k-l)\sum_{j=0}^{l-1}x^j
      \frac{{\rm d}^{k-l-1}}{{\rm d}x^{k-l-1}}f(x),$$
hence
$$\frac{{\rm d}^{k-l-1}}{{\rm d}x^{k-l-1}}f(x)=e^{(k-1)(x/1+x^2/2+\cdots+x^l/l)},$$
as required.
\end{proof}

Again, in view of Theorem~\ref{tp1} it is a challenge to find a
bijective proof for the number of $P_{1,l}^k$-avoiding
permutations in $S_n$ is the same number of $P_{a,l}^k$-avoiding
permutations in $S_n$, for all $1\leq a\leq k-l+1$.

\begin{example}\label{exp1}{\rm(see Mansour \cite{M})}
Theorem~\ref{tp1} for $l=1$ yields for $n\geq k$
    $$p_{a,1}^k(n)=(k-1)p_{a,1}^k(n-1),$$
with $p_{a,1}^k(k)=k!-(k-1)!$, hence
$p_{a,1}^k(n)=(k-2)!(k-1)^{n-k+2}$ for all $n\geq k-1$.\\
\end{example}

\begin{example}\label{exp2}{\rm(see Claesson and Mansour \cite{CM1} for $k=3$)}
Theorem~\ref{tp1} for $l=k-1$ yields for $n\geq k$
    $$\sum_{n\geq0}p_{a;k-1}^k(n)\frac{x^n}{n!}=e^{x/1+x^2/2+\dots+x^{k-1}/(k-1)}.$$
In particular, for $k=3$ yields the exponential generating
function for the number of permutations in $S_n(12\mn3, 21\mn3)$
{\rm(}or in $S_n(23\mn1,32\mn1)${\rm)} is given by $e^{x+x^2/2}$.
\end{example}

Again, as a remark, according to the main results
(Theorem~\ref{main11} and Theorem~\ref{tp1}), $C_{a,l}^k(x)$ and
$P_{a,l}^k$ does not depend on $a$; in other words,
$|S_n(C_{1,l}^k)|=|S_n(C_{a,l}^k)|$ and
$|S_n(P_{1,l}^k)|=|S_n(P_{a,l}^k)|$ for any $a$. We obtained this
fact as a consequence of lengthy computations. A natural question
would be to find a bijection between $S_n(C_{1,l}^k)$
(respectively, $S_n(P_{1,l}^k)$) and $S_n(C_{a,l}^k)$
(respectively, $S_n(P_{a,l}^k)$) that explains this phenomenon.
%========================================================================

\end{document}